\documentclass[12pt,leqno]{amsart}

\usepackage[dvips]{graphics}

\usepackage{amssymb}

\numberwithin{equation}{section}

\newtheorem{thm}{Theorem}[section]

\newtheorem{lem}[thm]{Lemma}
\newtheorem{cor}[thm]{Corollary}

\theoremstyle{definition}
\newtheorem{defn}{Definition}[section]
\theoremstyle{remark}
\newtheorem{rem}{Remark}[section]
\newtheorem{ex}[rem]{Example}

\newcommand{\tref}[1]{Theorem~\ref{#1}}
\newcommand{\cref}[1]{Corollary~\ref{#1}}

\newcommand{\lref}[1]{Lemma~\ref{#1}}

\newcommand{\R}{\mathbb{R}}

\newcommand{\Define}{\mathrel{\mathop:}=}

\begin{document}

%\tableofcontents

\pagebreak

%\bibliographystyle{alpha}

%\pagenumbering{roman}

\title{Spaces with many affine functions}

\author{Petra Hitzelberger and Alexander Lytchak }

%\twoauthors{Petra Hitzelberger}{Mathematisches Institut\\ Technische Universit\"at Darmstadt\\
%Schlossgartenstr.7, 64289 Darmstadt, Germany\\hitzelberger\@@mathematik.tu-darmstadt.de\\}{Alexander Lytchak}{Mathematisches Institut\\ Universit\"at Bonn\\
%Beringstr. 1, 53115 Bonn, Germany\\lytchak\@@math.uni-bonn.de\\}

%\footnotetext[1]{}


\begin{abstract}
We describe all metric spaces that have sufficiently many affine functions.
As an application we obtain a metric characterization of linear-convex 
subsets of Banach spaces.
\end{abstract}

\maketitle

\renewcommand{\theequation}{\arabic{section}.\arabic{equation}}
\pagenumbering{arabic}

\section{Introduction}
Given a geodesic metric space $X$ the existence of a non-constant
 affine function on $X$ seems
to impose strong  restriction on the geometry of the space.  Here a function 
$f:X\to \R$ is called affine if its restriction on each arclength parametrized
geodesic $\gamma$ in $X$ satisfies $f(\gamma (t))= at +b$ for some numbers
$a,b$ that may depend on $\gamma$. For instance,  if $X$ is a complete 
Riemannian manifold such a function happens to exist if and only if 
$X$ has a non-trivial Euclidean de Rham factor (\cite{inna}) and
the same  is true in some more general spaces (\cite{Mashiko},\cite{bishop}).  

 However, it is not 
clear what happens in the general situation. The expectation
 is  that the existence of  such a function forces the space $X$ to look
 like a Banach space in some "direction" defined by the affine function,
 compare the main result in \cite{LS} and examples  there. Thus it seems
 natural to expect that if $X$ has sufficiently many affine functions it should
 look very similar to some Banach space. Our result confirms this expectation:

\begin{thm} \label{main thm}
Let $X$ be a geodesic metric space. 
 Affine functions on $X$ separate points of $X$
if and only if $X$ is isometric to a convex subset of a normed vector
space with a strictly convex norm. 
\end{thm}
Here we say that affine functions separate points of $X$ if for each pair
of points $x, \bar x \in X$ there is an affine function $f:X \to \R$  with
$f(x) \neq f(\bar x)$.

It is possible to look on  our theorem from another point of view. 
In order to do this  we define a constant speed geodesic in a metric space 
to be a curve that has  a  
 constant speed and is globally minimizing between its endpoints. A map 
$F:X\to Y$ between geodesic metric spaces is called  affine if it sends each constant speed geodesic in $X$ to a constant speed geodesic
in $Y$.  We say that $F$ is an affine equivalence if $F$ is bijective and
$F$ and $F^{-1}$ are affine.

It is a natural question to which extent constant speed geodesic determine
the metric on a space, i.e. given a geodesic metric space $X$ what are all
metric spaces $Y$ that are affinely equivalent to the space $X$. If $X$ and $Y$
are complete Riemannian manifolds the answer   to this question  has been known
for a long time (\cite{inna}) and says that $Y$ has a 
de Rham decomposition that
consists of direct factors of $X$ stretched by constants.

\begin{rem}
If (still  in the realm of Riemannian geometry), one relaxes the condition
of affine equivalence to that of projective equivalence, i.e. if one
only requires the geodesics of $X$ and $Y$ to coincide
as subsets (unregarding the parametrisation) the question becomes much deeper.
In fact it is a very  old problem to which many  partial answers are known,
see, for instance, \cite{Matv} and the literature there.
\end{rem}

If (while looking for affinely equivalent spaces) 
we still ask $X$ to be a Riemannian manifold
but allow $Y$ to be arbitrary then the answer of \cite{inna}
must change. For example, the Euclidean space
is affine equivalent to each normed vector space  of the same 
dimension with a strictly
convex norm. As a consequence of our result we see that nothing more can happen
if $X$ is Euclidean. Namely  
the notion of an affine function does not change under affine equivalence, 
hence from \tref{main thm} we derive:

\begin{cor}
Let $C$ be  a convex subset of a Euclidean space $\R ^n$. If $Y$ is a geodesic 
metric space that
is affinely equivalent to $C$ then one can 
equip $\R ^n$ with a strictly convex norm
such that $C$ in this new norm is isometric to $Y$. 
\end{cor}

\begin{rem}
In  \cite{Ohta} it is shown that if $X$ is a Riemannian manifold 
and $Y$ is affine equivalent
to $X$, then  the metric on $Y$ is defined
by a continuous Finsler metric. We think that our arguments can  
help to describe the possible
Finsler metric precisely.
\end{rem}

The restriction to strictly convex norms in \tref{main thm}  seems artificial.
However, one should be careful when trying to drop this restriction.
In fact each metric space has an isometric embedding in some (non-strictly
convex) Banach space. On the other
hand $\R ^2$ with the maximum norm $||(x,y)||=\max \{ |x|,|y| \}$ has
too many geodesics to admit an affine function in the sense
of the definition above.  To avoid these difficulties one has to distinguish
good (i.e. linear) geodesics   from bad ones, thus we need
 the notion of a bicombing.

Let $X$ be a metric space. A bicombing $\Gamma$ on $X$ is an assignment 
to each pair $x,y$ of points in $X$ a geodesic $\gamma _{xy}$ connecting
$x$ and $y$, such that $\gamma _{yx}$ is the geodesic 
$\gamma _{xy}$ with the reverse orientation. 
Moreover, we require that for a point
$m$ on a geodesic $\gamma _{xy}$ the
 geodesic $\gamma  _{my}$ is part of the geodesic
$\gamma _{xy}$. We do not require that  the geodesics 
$\gamma _{xy}$ depend continuously on  $x$ and $y$.

 We say that a map $F:X \to Y$ between two spaces with 
bicombings denote by the same letter
 $\Gamma$  is $\Gamma$-affine if it for all $x,\bar x \in X$ there are some 
numbers $a,b$ such that 
 $f(\gamma _{x\bar x} (t))=\gamma _{f(x)f(\bar x)} (at +b)$ holds.

\begin{ex}
If $X$ is a uniquely geodesic metric space, i.e. if
between each two points of $X$ there  is precisely one geodesic,
then we have a unique  bicombing $\Gamma$ on $X$ and
 a $\Gamma$-affine function on $X$ is just an affine
 function in the sense of the old definition.
 \end{ex}

\begin{ex}
Each normed vector space has a natural bicombing that assigns to two points
the linear interval between 
them ($\gamma _{xy} (t) =x+ t \frac {y-x} {||y-x||})$.
All normed vector spaces will be considered with this bicombing $\Gamma$ 
in the sequel.
With this definition each linear map between two normed vector spaces 
is a $\Gamma$-affine
map.  On the other hand each $\Gamma$-affine map $F:V\to W$ between normed
vector spaces has the form $F(v)=F(0)+A(v)$ for some linear map $A$, i.e.
$F$ is affine in the usual sense of linear algebra. 
\end{ex}

 Using this notion we can now extend \tref{main thm}

\begin{thm} \label{exten}
Let $X$ be a space with a bicombing $\Gamma$. Then the following are 
equivalent:
\begin{enumerate}

\item There is a Banach space $B$ and an injective  isometric $\Gamma$-affine 
map $f:X\to B$;
\item The $\Gamma$-affine functions separate the points of $X$;
\item There is a Banach space $B$ and an injective  
$\Gamma$-affine  map $f:X\to B$.

\end{enumerate}
\end{thm}

The equivalence between (2) and (3) can be reformulated as a metric 
characterization
of linear convex subsets of Banach spaces:

\begin{cor}
Let $X$ be a metric space. Then $X$ is isometric to a linear convex subset
of some Banach space if and only if $X$ admits a bicombing $\Gamma$
such that the $\Gamma$-affine functions separate points of $X$.
\end{cor}

A short comment on the regularity of affine functions. It is well 
known that an affine functions on  a Riemannian manifold is smooth.
However, in general spaces affine functions do not need
to be even continuous, for instance, there are many non-continuous linear
functions on Hilbert spaces. In our theorem we did not make any restrictions 
on the affine functions, but a combination of \tref{exten} and the 
Hahn-Banach  theorem immediatly shows that for each geodesic metric space
 $X$ Lipschitz continuous 
affine function separate
points of $X$ if arbitrary affine functions separate points of $X$.

The paper is organized as follows. In Section \ref{easy} we show
that \tref{exten} implies \tref{main thm} and prove the straightforward
implications $(1)\to (2) \to (3)$ in \tref{exten}.
In Section \ref{reduct} we reduce the proof of the implication
$(3)\to (1)$ to the case where $X$ is an open convex subset of $\R^n$
with some Finsler metric on it.  Finally in Section \ref{vary} we discuss
first variation formulas in this  Finsler metric and prove that the 
Finsler structure is constant.

We would like to express our gratitude to Werner Ballmann for his 
encouragement and support. We are grateful to Linus Kramer 
for helpful comments.

\section{Preliminaries}

\subsection{Vector spaces}
Linear intervals in a vector space $V$ are curves $\gamma :[0,1] \to V$
of the form 
$\gamma (t) =t  v +(1- t) w$ 
for some $v,w\in V$. A subset $C$ of $V$
is called linearly convex if it  contains each linear  interval
$\gamma$ as above for all $v,w \in C$.  Let $C$ be a linearly convex subset
of $V$ that contains the origin $0$. 
Then the linear hull $H(C)$ of $C$ (i.e. the smallest linear subset of $V$
that contains $C$) is the set of all points $x$ that can be represented
as $x= \lambda (v- w)$ for some $\lambda >0$ and some  $v,w \in C$. 
 Let $C\subset V$ be linearly convex subset of $V$ that contains the origin
 and assume that $V$ coincides with the linear hull $H(C)$ of $C$. 
By the dimension of $C$ we  denote the dimension of the vector space 
$V=H(C)$.  If the dimension of $C$ is finite, it is well known that 
the set $O$ of inner points  
 of $C$ (with respect to
 the usual  Euclidean topology of $H(C)$) is convex, non-empty and that
the closure $\bar O$ of $O$ in $V$  contains $C$.

\subsection{Metric spaces} \label{Def und Notation}

By $d$ we will denote distances in metric spaces. 
A geodesic respectively a ray
in a metric 
space  $X$ will  denote an isometric
embedding $\gamma :I\to X$  of an interval respectively
 of a half-line into $X$.
Note that geodesics (if not otherwise stated) are parametrized 
by the arclength. 

 A metric space is called geodesic if each pair of its points is connected 
by a geodesic. It is called uniquely geodesic if this connecting geodesic 
is unique.

Given a ray $\gamma$ in a 
metric space $(X,d)$ the \emph{Busemann function of $\gamma$} is defined
by 
$b_{\gamma}(x)\Define \lim_{t \rightarrow \infty} 
(d(\gamma(t), x) - d(\gamma(0),\gamma(t)))$. 
The  Busemann function exists  and is a $1$-Lipschitz function on $X$.
Moreover $b_{\gamma} (x) \leq d(\gamma(t), x) - d(\gamma(0),\gamma(t))$
for all $t\geq 0$.

\subsection{Basics on normed vector spaces} \label{Busem}
Let $(V,\Vert \cdot \Vert )$ be a vector space with a norm and denote 
by $d$ the induced metric
on $V$. The norm is called strictly convex if for all linear
independent $v,w \in V$ one 
has $\Vert v+w \Vert < \Vert v \Vert + \Vert w \Vert$. The norm
$\Vert \cdot \Vert$ is strictly convex if and only if $(V,d)$ is uniquely
geodesic.

For each  non-zero vector $h$ in $V$  we denote by 
$\gamma ^h$ the 
ray $\gamma ^h (t) =t \frac h {\Vert h \Vert }$. The Busemann
functions $b_h := b_{\gamma ^h}$ of such linear rays  have
 the following properties, 
that are direct consequences of the definition.

\begin{enumerate}
\item  $b_h (v) =\lim _{t\to \infty} (||th -v|| -||th||) =
\lim _{t\to 0 ^+} (\frac {||h -tv|| -||h||} t)$;
\item $b_h(tv)=tb_h(v)$, for all  $t\geq 0$;
\item $b_h (v_1+v_2) \leq b_h (v_1) +b_h (v_2)$;
\item $b_h(v)=b_{-h}(-v)$.

\end{enumerate}

We see that the Busemann function $b_h$ is linear if and only if 
the equality  
$b_h (v)=-b_h(-v)$ holds for all $v \in V$. However, 
this equality
is equivalent to the  following well known one:

\begin{defn} \label{smdef}
A point $h\in V$ is called \emph{smooth in the Norm
$\Vert\cdot\Vert$ } if the following holds for all $v\in V$:
$$
\lim_{t \rightarrow 0^+}\frac{\Vert h -tv\Vert -\Vert h \Vert}{t} = -
\lim_{t \rightarrow 0^+}\frac{\Vert h +tv\Vert -\Vert h\Vert}{t}.
$$

\end{defn}

Thus the Busemann 
function $b_{h}$ is linear if and only if $h$ is a smooth vector
in the norm $\Vert \cdot \Vert$. If $V$ is a finite dimensional
vector space, then a vector $h$ is smooth in the norm $\Vert \cdot \Vert$
if and only if the Lipschitz map $\Vert \cdot \Vert :\R ^n \to \R$ has a linear
differential at the point $h$. 
Due to the theorem of Rademacher  we know that in a finite
dimensional normed vector space $V$ almost each  vector is smooth with
respect to the norm.

\section{Easy implications} \label{easy}

\subsection{Strictly convex case} First we are going to derive 
\tref{main thm} from \tref{exten}.  

Let $X$ be a space on which affine functions separate
the points. Let $x,\bar x$ be two arbitrary points in $X$. For each geodesic
$\gamma$ from $x$ to $\bar x$  each affine function $f$  has the value
$\frac {f(x)+f(\bar x)} 2 $ on the midpoint $m$ of $\gamma$. 
Thus all affine functions
have the same value on the midpoint of each geodesic between $x$ and $\bar x$. 
 From the 
separation assumption we deduce that $X$ is uniquely geodesic. 

 By \tref{exten} we know that $X$ is isometric to a linearly convex subset
of a Banach space $B$. We replace $X$ by this subset and may assume that
it contains the origin $0$. Denote by $V$ its linear hull in $B$ with the 
norm induced from the ambient Banach space $B$.  
Each element $v$ of $V$ has the form $\lambda (x-\bar x)$
for some $\lambda \geq 0$ and some $x,\bar x \in X$.

We claim that $V$ is strictly
convex. Assume the contrary and find linearly independent $v,w \in V$ with
$\Vert v+w \Vert =\Vert v \Vert + \Vert w \Vert $.   We can find a finite 
dimensional linearly convex subset $\tilde X$ of $X$ that contains $0$ and
 whose
linear hull $\tilde V$ contains $v$ and $w$. Due to finite dimensionality
 $\tilde X$
contains an open ball in $\tilde V$ and since $\tilde V$ is not strictly 
convex, the subset $\tilde X$ is not uniquely geodesic. But $\tilde X$ is a 
convex   subset of a uniquely geodesic space $X$ and we arrive at a 
contradiction.

The other direction is clear since linear functions are 
affine on a strictly convex normed vector space and separate its points
(compare the next subsection).

\subsection{Linear algebra} We prove  $(1)\to (2)$ of \tref{exten} here. 
Thus let   $(V, \Vert \cdot \Vert )$ be 
a normed vector space. Observe that 
each linear function on $V$ is $\Gamma$-affine with respect to the natural
bicombing. Since linear functions separate points of $V$ (by linear
algebra, you may also apply the Hahn-Banach theorem), 
we see that $\Gamma$-affine functions separate points of $V$.
We deduce that $\Gamma$-affine functions also separate points on each 
space $X$ that admits an injective  $\Gamma$-affine map $i:X\to V$.

\subsection{Evaluation} 
Assume now $(2)$ of \tref{exten}. Consider the vector space
$\Gamma Aff (X)$ of all $\Gamma$-affine functions on $X$ and denote
by $V$ its (algebraic) dual space $(\Gamma Aff (X) ) ^*$. 
Define the natural evaluation map $E:X\to V$ by $E(x) (f) = f(x)$
(compare \cite{LS}).
The map is well defined and by the separation property it is injective.  

 We claim that $E$  sends each $\Gamma$-geodesic onto a linearly reparametrized linear interval.

In fact, let a geodesic $\gamma _{xy}$ be given and reparametrize it to a constant speed geodesic $\bar \gamma _{xy}:[0,1] \to X$ such that 

$\bar \gamma _{xy} (0) =x, \bar \gamma _{xy} (1) =y$. For each $t\in [0,1]$ and each affine function $f:X\to \R$ we see $f(\bar \gamma _{xy} (t)) =(1-t) f(x) + tf(y)$. Hence for our map $E$ we deduce $E(\bar \gamma _{xy} (t) ) = (1-t) E(x) + tE(y)$.

Now we equip $V$ with an arbitrary norm. The linear intervals become linearly reparametrized $\Gamma$-geodesics for the natural bicombing $\Gamma$ on $V$. Hence the map $E:X\to V$ becomes  $\Gamma$-affine. Taking $B$ to be the completion of $V$ we obtain $(3)$ of \tref{exten}.

\section{Reduction} \label{reduct}

In this section we reduce the implication $(3)\to (1)$  of \tref{exten} to the case where $X$ is a Finsler manifold.

\subsection{Reformulation} 

Let $X$ be a space with  a bicombing $\Gamma$ that admits an injective $\Gamma$-affine map $i:X\to B$ into a Banach space $B$.  We forget the norm on $B$ and identify $X$ with its image $C$. We may assume that $0$ is contained in $C$ and replace $B$ by the linear hull $V$  of $C$ in $B$.  

We have to prove that we can define a norm on $V$ such that
the induced metric on $C$ coincides with the given one.
Thus the implication we are looking for is implied by the following:

\begin{lem} \label{reform}

Let $V$ be a vector space, $C$ a linear convex subset of $V$ that contains the origin $0$.
Let $d$ be a metric on $C$ such that the linear intervals contained in $C$ are constant speed geodesics with respect to the metric $d$. Then there is a norm $\Vert \cdot \Vert$ on $V$ such that the induced metric on $C$ coincides with $d$.
\end{lem}

\subsection{Reduction to finite dimension} 
 Since norms defined on linear subspaces can be extended
to norms on the whole space it is enough  to prove that such a norm exist on the linear hull $H(C) \subset V$.   Each element $v\in H(C)$ has the form $v= \lambda (x -\bar x)$  for some $x,\bar x \in C$ and some $\lambda \geq 0$.  Thus the norm we are looking for is unique, if it exists, and must be given by $\Vert v \Vert = \lambda  d(x,\bar x )$, where $v$ has a presentation as above.    

Thus we need to prove that this  is a well defined quantity (i.e. independent of the representation of $v$) and that it defines a norm on $H(C)$.

 Observe  that each linearly convex subset $\bar C$ of $C$ that contains $0$ and is considered with the metric induced by $d$ also satisfies the assumptions of \lref{reform}.
We  deduce that it is enough to prove \lref{reform} in the case where the dimension of $C$ is finite.   For instance, in order to see that $\Vert \cdot \Vert $ is well defined take a point $v\in H(C)$ that has two different presentations $v=\lambda (x- \bar x) $ and $v=\lambda _1 (x_1 -\bar x_1)$. Consider the linear convex hull $\bar C$ of the five points $x,\bar x, x_1, \bar x_1, 0$. Assuming that \lref{reform} is true for the finite dimensional  $\bar C$ we deduce that $\lambda  d(x, \bar x) = \lambda _1 d(x_1,\bar x_1)$, hence $\Vert \cdot \Vert $ is well  defined. The fact that $\Vert \cdot \Vert $ is a norm is  shown in the same way.

Thus it is enough to prove \lref{reform} in the case where $C$ has a finite dimension. 

\subsection{Reduction to open subsets} Let $C$ be as in \lref{reform} and assume that the dimension of $C$ is finite. Replacing $V$ by the linear hull $H(C)$ we may assume that $V=\R^ n$. If \lref{reform} is true for the set $O$ of inner  points of $C$ then by continuity it is also true for the closure $\bar O$ of $O$ in $\R ^n$ and therefore it is also true for $C\subset \bar O$. 

 Thus it is enough to prove \lref{reform} in the case where $V=\R^n$ and where $C$ is a convex open subset of $\R ^n$.

\subsection{Reduction to a Finsler metric}

Assume now that in  \lref{reform} we have $V=\R^ n$ and that $C$ is  open. We claim that the metric $d$ on $C$ is defined by a continuous Finsler structure on $C$. 

 In fact it is more or less a special case of Theorem B in \cite{Ohta}. 
We  shortly recall the arguments for the convenience of the reader, since the assumptions in \cite{Ohta} are slightly different. First of all one defines for each $x\in C$ and each $v\in \R ^n$  the quantity $|v| _x $ as $ d(x,x+\epsilon v ) /\epsilon$, for a sufficiently small positive real number $\epsilon$. Since $\gamma (t) =x +tv$ is a constant speed curve on the whole interval that is contained in the open set $C$, we deduce that the definition of $|v|_x$ does not depend on the choice of $\epsilon$, and in fact it is just the speed of the curve $\gamma$. 

Moreover, we immediately deduce that $|\lambda v|_x =|\lambda|\cdot |v|_x$. Since $d$ is a metric we have $|v|_x \geq 0$ with equality if and only if $v=0$.   

We  need to assure that the metric space $(C,d)$ is locally compact, i.e. that the identity $id:(C,d_{Eucl}) \to (C,d)$ is continuous. This follows directly, as soon as one knows that for  given $x\in C$ the function $|\cdot | _x$ is bounded on compact subsets of $\R ^n$ (or equivalently on the unit ball). The last claim  is the content of Proposition 2.1 in \cite{Ohta}. Fix now a small $\epsilon$ and take a convergent sequence $x_i \to x$ and a convergent sequence $v_i \to v$. We observe that $x_i+\epsilon v _i $ converge to $x+\epsilon v$, hence by definition $|v_i|_{x_i} $ must converge to $|v|_x$. Therefore the function $|v|_x$ is continuous in $v$ and $x$. From this we can derive that $|\cdot |_x$ is a norm. Namely, we have 
$d(x,x+\epsilon (v +w)) \leq d(x,x+\epsilon v) +
d(x+\epsilon v,x+\epsilon v + \epsilon w) =
\epsilon |v|_x + \epsilon |w| _{x+\epsilon v}$. Dividing by $\epsilon$ we get in the limit $|v+w|_x \leq |v|_x + |w|_x$.

Thus $|\cdot |_x$ is a continuous Finsler structure and defines a metric  $\tilde d$ on $C$.  The conclusion that $\tilde d$ coincides with $d$ is a direct consequence of the fact that linear intervals have the same length with respect to the metric $d$ and $\tilde d$.

Thus we have shown that the metric $d$ is given by a continuous Finsler structure $|\cdot |_x$ on $C$.

\section{Final step} \label{vary}

We are left with the following problem. Let $O\subset \R ^n$ be an open linear convex subset. Let $|\cdot |_x$ be a continuous Finsler structure on $O$, such that for the induced metric $d$ the  linear intervals are constant speed geodesics. We need to show that there is a norm $\Vert \cdot \Vert $ on $\R ^n$ such that the metric $\tilde d$ on $O$ induced by this norm coincides with our metric $d$. This conclusion follows as soon as we know that $|\cdot |_x =|\cdot | _y$ holds for all $x,y\in O$.   The rest of this section is devoted to the proof of the last statement.

 First of all  linear intervals have constant speed, hence for each $x\in O$, each $v\in \R ^n$ and each $t\in \R$ such that $x+tv$ is still in $O$ we have $|v| _x =|v|_{x+tv}$. Since linear intervals are geodesics we conclude that $d(x,x+tv) =|tv| _x$, i.e. for all $x,y \in O$ we have $d(x,y) =|x-y|_x$. (To avoid confusion: what we need is the much stronger statement $d(y,z)=|y-z| _x$ for all $x,y,z \in O$!).   

Let us fix an arbitrary $v\in \R^n$. It is enough to show that $|h|_x=|h|_{x+v}$ for almost all $h \in \R ^n$ and all $x,x+v \in O$. Here and below we consider $ \R^n, O$ the tangent bundle $TO= O\times \R^n$ and linear intervals equipped with natural Lebesgue measures.
By abuse of notation we will say that $h\in \R ^n$ is smooth at a point $x\in O$ if $h$ is a smooth vector of the norm $|\cdot | _x$ (compare Definition \ref{smdef}). For each $x\in O$ almost all $h \in \R^n$ are smooth at $x$. 
Applying Fubini's theorem twice we find a set $S$ of full measure  in $\R ^n$ with the following property.  For each vector $h\in S$ the set of all $x\in O$ at which $h$ is smooth  is of full measure in $O$.

It is enough to show that for each $h\in S$ the equality  $|h|_x=|h|_{x+v}$ holds, for all $x,x+v \in O$. Let us fix a vector $h\in S$.  Another application of Fubini's theorem tells us, that for almost all $x\in O$ the vector $h$ is smooth at almost each point of the  linear interval $\eta (t) =x+tv$ (contained in $O$). By continuity of the Finsler structure it is enough to show that $f(t) =|h|_{x+tv}$ does not depend on $t$ if $h$ is smooth at almost each point of $\eta$. The question is local in  $t$, hence it is enough to show that $f(t)$ is constant for all small $t$. We may assume (replacing $h$ by its multiple) that $x+h$ is contained in $O$. Then $f(t)=d(x+tv,x+tv+h)$, in particular, the function $f$ is locally Lipschitz. Hence it is enough to show that $f'(t) =0$ for almost all $t$. Therefore our result is a consequence of the following lemma.

\begin{lem} \label{final}
Let $x$ be a point in $O$, $h,v\in \R^n$. Assume that $h$ is a smooth vector in the norm $|\cdot|_x$. Then for the function $f(t)= |h|_{x+tv}$ we have $f'(0) =0$.
\end{lem}

\begin{proof}
Replacing $h$ by its multiple we may assume that $y:=x+h$ is contained in $O$.  Our claim is $|h|_{x+tv} =|h|_x +o(t)$, where $o(t)$ is some function with $\lim _{t\to 0} \frac {o(t)} t =0$. Since $|h|_{x+tv} =d(x+tv,y+tv)$ we just have to find right first variation formulas in the Finsler manifold $O$. We refer to Section 9 in \cite{alex1} for a general discussion of first variation formulas.
We need some notations. For a vector $w \in \R ^n$ a point $z\in O$ we will denote by $b_w ^z :\R^n \to \R$ the Busemann function of the direction $w$ in the norm $|\cdot |_z$, i.e. $b_w ^z (v) = \lim _{t\to \infty} (|tw-v|_z -|tw|_z)$. We refer to Subsection~\ref{Busem} for basic properties of such functions.

First we are going to show that $b_h ^x =b_h ^y$. In order to do so we have to study the first derivative of the function $\tilde f (t)= d(y,x+tv)$. Recall that $|h|_x =|h|_y$.

Due to  $d(y,x+tv)=|y-x-tv|_y=|h-tv|_y $  and properties of the Busemann functions (see Subsection \ref{Busem}) we have 
$$
\text{(1)\hspace{2ex}} d(y,x+tv)=  |h|_y +t b^y _{h} (v) + o(t). 
$$
 On the other hand we can estimate $d(y,x+tv)$ by working  
in a small neighborhood of $x$. Fix $s \gg 1$. For $0<t \ll \frac 1 s$ we have $d(y,x+tv) \leq d(y,x+tsh) + d(x+tsh,x+tv)$.  Since the metric on $O$ is given by a continuous Finsler structure, we have $d(x+tsh,x+tv) =t |sh-v|_x + o(t)$ (i.e. $(\R ^n ,|\cdot |_x)$ is the tangent space of $O$ at $x$ in the sense of metric geometry, \cite{alex1}). We obtain $d(y,x+tv) \leq |h|_x (1-ts) + t |sh-v|_x  + o(t) $.
If $s$ is large enough $|sh-v|_x $ is very close to $s +b^x _{h}(v) $. Letting $s$ go to infinity we conclude, that 
$$
\text{(2)\hspace{2ex}} d(y,x+tv) \leq |h|_x + t b^x _{h^+} (v) + o(t)  \text{\hspace{2ex} for $t>0$.}
$$
Consider 
$d(y,x+tv)+ d(x+tv, x-tsh) \geq  d(y,x-tsh) =|h|_x(1-ts)$
for a fixed big $s$ and let $t$ go to $0^+$.  Then as above we obtain  
$$
\text{(3)\hspace{2ex}} d(x+tv,y) \geq |h|_x - t b ^x _{-h} (v)  + o(t)  \text{\hspace{2ex} for $t>0$.} 
$$
By our assumption vector $h$ is smooth at $x$, thus $b^x _{-h} (v) =- b ^x _{h} (v)$. Combining $(1),(2),(3)$  above we deduce that $b^x_{h} (v) =b^y _{h} (v)$. Since the above calculations work for all $v\in \R^ n$  we get $ b^x_{h}  =b^y _{h} $. In particular, $h$ is also smooth in $|\cdot |_y$ and we have  $ b^x_{-h}  =b^y _{-h}$. Now we can finish the proof of \lref{final}. Arguing as in the proof of (2) and (3) above. We take a large $s$ and estimate the distance $d(x+tv,y+tv)$ from above (for all small  $0<t \ll \frac 1 s$) by 
\begin{align*}
  d(x+tv,y+tv) \leq &  d(x+tv, x+tsh)+ d(x+tsh, y-tsh)\\
                    & +d(y-tsh , y+tv) 
\end{align*}
Due to the continuity of the Finsler structure we have
$$ d(x+tv,x+tsh) =t|v-sh|_x +o(t) $$
$$ \text{and}\hspace{2ex} d(y-tsh , y+tv) = t |v +sh| _y  +o (t). $$
Thus we conclude
$$
d(x+tv,y+tv)\leq  t|v-sh| _x +|h| (1 -2ts)  +t |v+sh|_y + o(t). 
$$
If $s$ goes to infinity then $|v-sh|_x -s$ converges to 
$ b^x _{h} (v)$
and  $|v+sh|_y -s$ converges  to  $b ^y _{-h} (v)$. Thus we arrive at:
$$
\text{(4)}\hspace{2ex} d(x+tv,y+tv) \leq |h| + t b^x _{h} (v)  +t b ^y _{-h} (v) +o(t).
$$
  Now we  take a large $s$ and estimate $d(x+tv,y+tv)$ from below by 
\begin{align*}
d(x+tv,y+tv)\geq &  d(x-tsh,y+tsh) -d(x-tsh, x+tv)\\
                 & - d(y+tv, y+tsh)
\end{align*}
As above we derive from it
$$
\text{(5)}\hspace{2ex}  d(x+tv,y+tv)  \geq |h| - t b^x _{-h} (v) - t b^y _{h} (v) + o(t).
$$
Since we know that $b^x_{h} =b^y_{h}=-b^y_{-h}=-b^x_{-h}$ we can combine  $(4)$ and $(5)$ and obtain $d(x+tv,y+tv) =d(x,y) +o(t)$. Thus we have shown $f' (0) =0$.
\end{proof}

This finishes the proof of \tref{exten}.

\bibliographystyle{alpha}
\bibliography{meine}

\begin{tabular}{ll}
&\\
\em{Petra Hitzelberger} & \em{ Alexander Lytchak} \\
{\small Mathematisches Institut}  &{\small  Mathematisches Institut} \\ 
{\small TU Darmstadt}     &{\small Universit\"at Bonn }\\
{\small Schlossgartenstr.7, 64289 Darmstadt} & 
{\small Beringstr. 1, 53115 Bonn}\\
{\small Germany }  &{\small  Germany }\\
{\tiny hitzelberger\@@mathematik.tu-darmstadt.de} & {\tiny lytchak\@@math.uni-bonn.de}\\ 
\end{tabular}

\end{document}